\documentclass[reqno,a4paper,11pt]{amsart}
\usepackage{amssymb}
\usepackage[dvips, lmargin=4cm, rmargin=4cm, tmargin=4cm, bmargin=4cm]{geometry}
\usepackage[colorlinks=true, citecolor=red, linkcolor=blue, urlcolor=blue]{hyperref}
\usepackage{cite}
\usepackage{xfrac}
\usepackage{graphicx}
\usepackage{hyperref}
\usepackage{amsmath}
\usepackage[T1]{fontenc}
\usepackage{ulem}
 \usepackage{array,multirow,makecell}
 \setcellgapes{1pt}
 \makegapedcells
 \newcolumntype{R}[1]{>{\raggedleft\arraybackslash }b{#1}}
 \newcolumntype{L}[1]{>{\raggedright\arraybackslash }b{#1}}
 \newcolumntype{C}[1]{>{\centering\arraybackslash }b{#1}}
\usepackage{enumitem}
\usepackage{dsfont}
\everymath{\displaystyle}
\newtheorem{theorem}{Theorem}[section]

\newtheorem{remark}{Remark}[section]
\newtheorem{lemma}{Lemma}[section]
\newtheorem{proposition}{Proposition}[section]
\newtheorem{cor}{Corollary}[section]

\numberwithin{equation}{section}
\newtheorem{notation}{Notation}[section]

\begin{document}
\title[Fractional anisotropic problem with variable exponent]{Anisotropic fractional Sobolev spaces with variable exponent and application to nonlocal problems}
\author[E. Azroul, A. BARBARA, N. Kamali and  M. Shimi]
{E. Azroul$^1$, A. Barbara$^2$, N. Kamali$^3$ and M. Shimi$^4$ }
\address{Elhoussine Azroul, Abdelkrim BARBARA and Nezha KAMALI\newline
 Sidi Mohamed Ben Abdellah
 University,
 Faculty of Sciences Dhar el Mahraz, Laboratory of Mathematical Analysis and Applications, Fez, Morocco}
 
 \email{$^1$elhoussine.azroul@usmba.ac.ma}
  \email{$^2$abdelkrim.barbara@usmba.ac.ma} 
 \email{$^3$nezha.kamali@usmba.ac.ma}
\address{Mohammed Shimi\newline
 Sidi Mohamed Ben Abdellah
 University,
ENS of Fez, Laboratory of Mathematical Analysis and Applications, Fez, Morocco}
\email{$^4$mohammed.shimi2@usmba.ac.ma}
\date{Month, Day, Year}
\begin{abstract}
The main goal of this paper is to introduce a new fractional anisotropic Sobolev space with variable exponent where the basic qualitative properties (completeness, separability, reflexivity, ...) are established, including the continuous and compact embedding results. Moreover, some functional proprieties of anisotropic fractional $\vec{p}(.,.)$-Laplacian operator are proved.  As an application, we use the  mountain pass theorem and Ekeland's variational principle to ensure the existence of a weak solution for a nonlocal anisotropic problem with variable exponent.  
\end{abstract}
\maketitle
\tableofcontents
\section{Introduction}
The theory of Sobolev spaces, $W^{k,p}(\Omega)$ consists of functions defined on a domain $\Omega$ that have weak derivatives up to order $k$ and are $L^p$-integrable, where $p$ is a real number greater than or equal to 1,  has been originated by the mathematician S.L. Sobolev in \cite{Sobolev} around 1938. These  spaces are a fundamental concept in functional analysis and partial differential equations (PDEs). They provide a framework for studying PDEs and to prove the  regularity, existence, and uniqueness of solutions to a wide range of differential equations because they allow to extend the notion of a solution to equations where traditional smooth solutions may not be well-defined. 

The variability in regularity may arise in various scientific and engineering applications, then one of the key motivations for introducing variable exponents $p(x)$ is to model situations where different parts of a domain exhibit distinct degrees of regularity. Hence, the introduction of variable exponent Sobolev spaces allows for a more flexible and nuanced analysis of problems where the regularity of solutions is not uniform across the entire domain, providing a valuable extension of classical Sobolev space, for a deep comprehension of the  theory of these spaces, we refer to  the monograph  \cite{diening} and the references \cite{Fan, Kovacik}. Naturally, problems involving the $p(x)$-Laplace operator
$$
\Delta_{p(x)} (u)=\operatorname{div}\left(|\nabla u|^{p(x)-2} \nabla u\right)
$$
were extensively investigated. In parallel, due to the development of various real applications, including physics, materials science, image processing, and computational mathematics the term "anisotropic exponent" was introduced. The inception of this concept of different directions was by Benedek and Panzone \cite{benedek}, where they investigated the mixed Lebesgue spaces. Originally, this birth date back to \cite{hormander} in the estimation for translation invariant operator in Lebesgue spaces, later on, fundamental studies in this direction are done (see for example \cite{igari, madych, sjodin}). In that context, in 1969 the notion of anisotropic Sobolev spaces was introduced by Troisi \cite{troisi}, $W^{k, \vec{p}}$ (where $\vec{p}$ is a constant vector, $\vec{p}=\left(p_1, \ldots, p_N\right)$), there are also many scholars \cite{kolodi, mihalescu, bedahmane} dealing with problems involving the $\vec{p}$-Laplace operator,
$$
\Delta_{\vec{p}}(u)=\sum_{i=1}^N \partial_{x_i}\left(\left|\partial_{x_i} u\right|^{p_i-2} \partial_{x_i} u\right).
$$ 
 Furthermore, a new theory riveted attention when it introduced the anisotropic space with variable exponent, $W^{1, \vec{p}(\cdot)}$ (where $\vec{p}(\cdot)=\left(p_1(\cdot), \ldots, p_N(\cdot)\right)$ is a vector with variable components) see \cite{fan1}. As a consequence, a new operator takes its place in the   literature, explicitly
$$
\Delta_{\vec{p}(x)}(u)=\sum_{i=1}^N \partial_{x_i}\left(\left|\partial_{x_i} u\right|^{p_i(x)-2} \partial_{x_i} u\right) ;
$$
see also \cite{mihalescu,mihalescu2} and the references therein. This operator is referred as the anisotropic variable exponent $\vec{p}(\cdot)$-Laplace operator and it is intimately connected with $\Delta_{{p}(x)}$ and $\Delta_{\vec{p}}$. 
 
During the same period, the theory of fractional calculus has piqued the interest of many scientists and engineers since fractional differential equations, which involve fractional derivatives, are used to describe processes that exhibit memory effects and long-range dependencies and due to their application in widespread fields such as image processing \cite{yang}, electro-chemistry \cite{Oldham}, and electromagnetic \cite{engheta}, etc. In particular,  the fractional Laplacian operator is a powerful mathematical tool that extends the classical Laplacian to non-integer orders, It is a typical example of a non-local operator, meaning its action at a point depends on the values of the function over the entire domain. This non-locality is crucial for capturing behaviors that extend over long distances or are influenced by distant points. Thence, nonlocal equations associated with this type of operator became one of the prevalent fields of investigation both for pure mathematical research and in view of concrete
real-world applications, we refer to \cite{dinezza} for a full introduction to the study of this operator and the fractional Sobolev spaces related to them. Let us fixing $s$ in $(0,1)$, for any $p \in$ $[1,+\infty)$, the space $W^{s, p}(\Omega)$ is defined as follows
$$
W^{s, p}(\Omega):=\left\{u \in L^p(\Omega): \frac{|u(x)-u(y)|}{|x-y|^{\frac{n}{p}+s}} \in L^p(\Omega \times \Omega)\right\} .
$$
These spaces are intermediary Banach spaces between $L^p(\Omega)$ and $W^{1, p}(\Omega)$ and they associated with the fractional $p$-Laplacian operator $(-\Delta)^s_p$.\\
 On the one hand, the authors in \cite{Xu},  elongated these spaces to the anisotropic fractional Sobolev spaces and they proved that
such spaces embed to mixed-norm Lebesgue spaces. Moreover, they also give some embedding theorems between these spaces. Otherwise, the space $W^{s, p}$ was extended to the variable exponent case where many researchers have made remarkable achievements in the study of the fractional Sobolev spaces with variable exponent $W^{s, p(x,y)}$ and the related fractional partial differential equations with variable growth involving the fractional $p(x,.)$-Laplacian $(-\Delta)^s_{p(x,.)}$ (see for example \cite{ayazoglu,azroul2,azroul1,azroul3,bahrouni,ho,kaufmann,Liu,Zuo}).\\

Motivated by the above works,  the challenge now is whether we can extend this concept of anisotropy to the fractional Sobolev spaces with variable exponent.  To this aim, we want to define a new kind of functional framework that covers all the aforementioned spaces.  These new functional spaces provide a framework for studying functions in different directions; to be more precise, they will be used to deal with functions with varying regularity in different directions.\\

Let $\Omega$ be a smooth bounded open domain  in $\mathbb{R}^N$. We consider a continuous vectorial function $\vec{p}:\left( \overline{\Omega} \times \overline{\Omega} \right)^N \longrightarrow \mathbb{R}^N$   which verifies, for all $ i=1,...,N$, 
\begin{equation} \label{p1}
 {\small 1<p_i^{-}=\min _{(x, y) \in \overline{\Omega} \times \overline{\Omega}} p_i(x, y) \leqslant p_i(x, y) \leqslant p_i^{+}=\max _{(x, y) \in \overline{\Omega} \times \overline{\Omega}} p_i(x, y)<+\infty},  
 \end{equation}
 and 
 \begin{equation} \label{p2}
 p_i \text{ is symmetric, i.e, } p_i(x, y)=p_i(y, x) \quad \text{ for all } (y, x) \in\overline{\Omega} \times \overline{\Omega} .
 \end{equation}
 To make the work easier for the reader, we have included a list of notations to avoid confusion between variables. 
 \begin{notation}
 for all $(x,y) \in \overline{\Omega}\times \overline{\Omega}$, we denote\\
  \begin{minipage}{8cm}
\begin{tabular}{l}
$\bullet$ $\vec{p}(x,y)= \left\lbrace  p_1(x,y), p_2(x,y),...,p_N(x,y)\right\rbrace ,$\\
 $\bullet$ $\overline{p}_i(x)= p_i(x,x)$,\\
 $\bullet$ $p_M(x)= \underset{1\leqslant i \leqslant N}{\max} \overline{p}_i(x)$,\\
  $\bullet$ $\vec{P}_{+}=\left(p_1^{+}, \ldots, p_N^{+}\right),$\\
  $\bullet$ $ \vec{P}_{-}=\left(p_1^{-}, \ldots, p_N^{-}\right)
    $,
\end{tabular}  
  
 \end{minipage}
  \begin{minipage}{5.3cm}
 \begin{tabular}{l}
 $\bullet$ $
  P_{+}^{+}=\max \left\{p_1^{+}, \ldots, p_N^{+}\right\},$  \\
  $\bullet$  $P_{+}^{-}=\min \left\{p_1^{+}, \ldots, p_N^{+}\right\},$  \\
  $\bullet$  $ P_{-}^{+}=\max \left\{p_1^{-}, \ldots, p_N^{-}\right\}, $ \\
   $\bullet$  $P_{-}^{-}=\min \left\{p_1^{-}, \ldots, p_N^{-}\right\} .
    $    
 \end{tabular}  
   
  \end{minipage}  

 \end{notation}
Now, for $s\in(0,1)$, we introduce the fractional anisotropic Sobolev space with variable exponent as follows
{\small  $$
W=W^{s, p_{{\scriptsize M}}(x),\vec{p}(x,y)}(\Omega)
 =\left\{ \begin{aligned}
 w \in & L^{p_M(x)}(\Omega) :  \iint_{\Omega\times\Omega}\tfrac{\left| w(x)-w(y)\right|^{p_i(x,y)} }{\nu_i^{p_i(x,y)}\left| x-y \right| ^{sp_i(x,y)+N} } dx dy <+\infty, \\
 & \quad \text{for some } \nu_i>0, \text{ and for all }i=1,...,N
 \end{aligned} \right\}.
 $$}
The space $W$ can be equipped with this norm  
\begin{equation}
\|w\|_{s,\vec{p}(x,y)} =\left[ w\right]_{s,\vec{p}(x,y)}+ \|w\|_{p_{M}(x)}, 
\label{norm}
\end{equation}
where $\left[ w\right]_{s,\vec{p}(x,y)}$ is defined as the sum of the Gagliardo seminorm with variable exponent and it is given by 
\begin{equation} \left[ w\right]_{s,\vec{p}(x,y)}= \sum_{i=1}^{N} \inf \left\{\nu_{i}>0: \iint_{\Omega \times \Omega} \frac{|w(x)-w(y)|^{p_{i}(x, y)}}{\nu_{i}^{p_{i}(x, y)}|x-y|^{sp_{i}(x, y)+N}}\leqslant1\right\} . 
\end{equation}
These spaces are intermediary spaces between $L^{\vec{p}(x)}(\Omega)$ and $W^{1, \vec{p}(x)}(\Omega)$.\\ 
It is worth mentioning that in the literature, there is one paper \cite{azroul4} in which an anisotropic fractional $\left( p_1\left( x,.\right) ,p_2\left( x,.\right)\right)$-Kirchhoff type problems investigated. However, to the best of our knowledge, the fractional anisotropic Sobolev spaces with variable exponent have not yet been studied systematically up to now. The main goal of this article is to present some basic properties of these spaces and the fractional $\vec{p}(x,.)$-Laplacian operator
\begin{equation}
{\small  \left[ \left( -\Delta\right)^{s}_{\vec{p}(x,.)} u\right] (x)=p . v \sum_{i=1}^{N} \int_{\mathbb{R}^N} \tfrac{|u(x)-u(y)|^{p_{i}(x, y)-2}(u(x)-u(y))}{|x-y|^{N+s p_{i}(x, y)}} d y, \quad \text{for all } x\in \mathbb{R}^N},
\end{equation}
where p.v. is a widely used shorthand for principal value,
 so as to provide a basis for the further study of the nonlocal equations involving this operator.\\
 
Note that when choosing $p_1(\cdot,\cdot)=\cdots=p_N(\cdot,\cdot)=p(\cdot,\cdot)$ we obtain an operator with similar properties to the fractional $p(x,\cdot)$-Laplacian operator $(-\Delta)^s_{{p}(x,.)}$. Furthermore $\left( -\Delta\right)^{s}_{\vec{p}(x,.)}$ is the fractional version of the $\vec{p}(x)$-Laplacian operator, while when choosing $p_1, \ldots, p_N$ to be constant functions we arrive at the anisotropic fractional $\vec{p}$-Laplacian operator $(-\Delta)^s_{\vec{p}}$.\\

A fundamental result for Sobolev spaces is the embedding theorem. One of the main goals of this paper is to establish an embedding theorem that consists that the anisotropic fractional Sobolev spaces with variable exponent are continuously embedded in a generalized Lebesgue space. This result is stated as follows  
\begin{theorem}\label{th1}
Consider a smooth bounded domain $\Omega$ in $\mathbb{R}^N$ $(N\geqslant 2)$, and let $\vec{p}:\overline{\Omega}\times \overline{\Omega} \longrightarrow (1,+\infty)$ be a vector function satisfying conditions \eqref{p1} and \eqref{p2} with $sp_i(x,y)<N$, $ i=1...N,$  for all $(x,y)\in \overline{\Omega}\times \overline{\Omega}$. Moreover, assume that $q: \overline{\Omega}\longrightarrow   (1,+\infty)$ is a continuous variable exponent such that 
$$ 1 < q(x)< p_i^{s,*}(x)=\frac{p_i(x,x)} {1-\frac{sp_i(x,x)}{N} }, \quad i=1, \ldots, N, \text { for all } x \in \overline{\Omega}.     $$
Therefore, there exists a constant $C$ that depends on $N$, $s$, $\vec{p}$, $q$, and $\Omega$ such that 
\begin{equation}
\|u\|_{L^{q(x)}(\Omega)} \leqslant C\|u\|_{s,\vec{p}(x,y)}.
\label{ine1}
\end{equation}
Hence, the embedding of $W^{s, p_M(.),\vec{p}(.,.)}\left( \Omega\right) $ into $L^{q (x)}(\Omega)$ is continuous. Furthermore, this embedding is compact.
\end{theorem}

The remainder of this paper is structured as follows, Section \ref{S2} recalls the necessary theoretical background of Lebesgue and fractional Sobolev spaces with variable exponent, In section \ref{S3}, we discuss the crucial properties of the variable exponent anisotropic fractional Sobolev space   such as the completeness, separability, reflexivity. Besides, we give a detailed proof of the  continuous and compact  embedding result. In addition, we will adjust the definition of these spaces in order to encode correctly the Dirichlet boundary datum in nonlocal problems.
 In section \ref{S4}, we establish some basic   properties of the fractional $\vec{p}(x,.)$-Laplacian functionals.
 Section \ref{S5} is devoted to proving the existence of weak solution of an application associated with our abstract setting involving the fractional $\vec{p}(x,.)$-Laplacian operator  using some variational methods.
\section{Preliminary results}\label{S2}
This section is devoted to briefly reviewing some basic properties concerning the  Lebesgue and  fractional Sobolev spaces with variable exponents.
 \subsection{Variable exponent Lebesgue spaces}
Consider a bounded subset $\Omega$ of $\mathbb{R}^N$ $(N\geqslant 2)$, and let $ C^{+} 
(\overline{\Omega})$ defined as follows
 $$ C^{+} 
(\overline{\Omega})= \{ p \in C(\overline{\Omega}): 1< p(x) ~~ \text{  for all } x \in \overline{\Omega}  \}.$$
For any $p \in C^{+} (\overline{\Omega})$, we denote 
$$ p^{+}=\underset{x\in \overline{\Omega}}{\sup}~p(x) ~~ \text{  and   }  ~~ p^{-}=\underset{x\in \overline{\Omega}}{\inf}~p(x), \quad \text{where} \quad 1<p^-\leqslant p(x) \leqslant p^+<+\infty.$$
Then, for any $p \in C^{+} (\overline{\Omega})$, the Lebesgue space with variable exponent is characterized by the following definition:
$$ L^{p(x)}(\Omega)= \left\{w  \in \mathcal{M}(\Omega,\mathbb{R})  \text{ such that: } \rho_{p(.)}(w):=\int_{\Omega}|w(x)|^{p(x)} d x<\infty\right\}. $$
The Lebesgue space with variable exponent endowed with this norm $$\|w\|_{p(x)}=\|w\|_{L^{p(x)}(\Omega)}=\inf \left\{\nu>0 ; \int_{\Omega}\left|\frac{w(x)}{\nu}\right|^{p(x)} \mathrm{d} x \leqslant 1\right\}
  ~~~~~$$ is a separable, reflexive, and Banach space.
   \begin{proposition}[\cite{Fan}]\label{pro1}
    For all $w \in L^{p(x)}(\Omega)$, we define the modular  $$\rho_{p(x)}(w)=\int_{\Omega}\left|{w(x)}\right|^{p(x)} \mathrm{d} x.$$ Then, we have the following relations between the norm and modular.
      \begin{enumerate}
               	\item $\|w\|_{L^{p(x)}(\Omega)}=1$ $($resp. $<1,>1) \Leftrightarrow \rho_{p(x)}(w)=1$ $($resp. $<1,>1)$,\\
         \item $\|w\|_{L^{p(x)}(\Omega)}>1 \Rightarrow\|w\|_{L^{p(x)}(\Omega)}^{p-} \leqslant \rho_{p(x)}(w) \leqslant\|w\|_{L^{p(x)}(\Omega)}^{p+}$, \\
          \item $\|w\|_{L^{p(x)}(\Omega)}<1 \Rightarrow\|w\|_{L^{p(x)}(\Omega)}^{p+} \leqslant \rho_{p(x)}(w) \leqslant\|w\|_{L^{p(x)}(\Omega)}^{p-}$, \\
    	\item $\lim _{k \rightarrow+\infty}\left\|w_{k}-w\right\|_{L^{p(x)}(\Omega)}=0 \Leftrightarrow \lim _{k \rightarrow+\infty} \rho_{p(x)}\left(w_{k}-w\right)=0 .$
                	
                	\end{enumerate}
    \end{proposition}
\begin{lemma}[\cite{Fan}]\label{Holder}
For all $v \in L^{p(x)}(\Omega)$ and $w \in L^{\hat{p}(x)}(\Omega)$, with $\tfrac{1}{p(x)}+ \tfrac{1}{\hat{p}(x)}=1 $, then
           \begin{equation}\label{Hlr}
           \left|\int_{\Omega}  v w \mathrm{~d} x\right| \leqslant\left(\tfrac{1}{p^{-}}+\tfrac{1}{(\hat{p})^-}\right)\|v\|_{L^{p(x)}(\Omega)} \|w\|_{L^{\hat{p}(x)}(\Omega)} \leqslant 2\|v\|_{L^{p(x)}(\Omega)} \|w\|_{L^{\hat{p}(x)}(\Omega)},
           \end{equation}
it is the so-called H\"{o}lder's inequality.          
\end{lemma}   
 An immediate consequence of H\"{o}lder's inequality is the following.
 \begin{cor}[\cite{Fan}]\label{H2.1}
 If $r(.), q(.)\in C_+(\overline{\Omega})$, define $p(.)\in C_+(\overline{\Omega})$ by $$\frac{1}{p(x)}=\frac{1}{q(x)}+\frac{1}{r(x)}.$$
 Then there exists a positive constant $C$ such that for all $u\in L^{q(x)}(\Omega)$ and $v\in L^{r(x)}(\Omega)$ we have that $uv\in L^{p(x)}(\Omega)$ and
 $$\|u v \|_{L^{p(x)}(\Omega)}\leqslant C\|u\|_{L^{q(x)}(\Omega)}\|v\|_{L^{r(x)}(\Omega)}.$$ 
 \end{cor} 
  \subsection{Fractional Sobolev space with variable exponent}
 For any $s \in (0,1)$ and $p \in C^{+} (\overline{\Omega} \times \overline{\Omega}) $ that verifies conditions \eqref{p1} and \eqref{p2} (i.e. for $i=1$), the fractional Sobolev space with variable exponent is defined as follows 
 $$\begin{aligned}
 Y&=W^{s, p(x, y)}(\Omega) \\
  &=\left\{w \in L^{\bar{p}(x)}(\Omega): \int_{\Omega \times \Omega} \tfrac{|w(x)-w(y)|^{p(x, y)}}{\nu^{p(x, y)}|x-y|^{s p(x, y)+N}} d x d y<+\infty, \text { for some } \nu>0\right\}.
 \end{aligned}$$
We equip this space with the following norm 
 $$\|w\|_{Y}=\|w\|_{s, p(x, y)}=\|w\|_{L^{\bar{p}(x)}(\Omega)}+[w]_{s,p(x,y)},$$
 with \begin{equation}
 \left[ w\right]_{s,p(x,y)}=  \inf \left\{\nu >0: \int_{\Omega \times \Omega} \frac{|w(x)-w(y)|^{p(x, y)}}{\nu^{p(x, y)}|x-y|^{sp(x, y)+N}}\leqslant 1\right\} .
 \label{gagliardo}
 \end{equation} 
 Therefore, the space $(Y,\|w\|_{Y})$ transformed to a separable and reflexive Banach space (we refer the reader to \cite{azroul1,kaufmann,bahrouni}).\\
 We denote by $Y_0=W^{s, p(x, y)}_0(\Omega)$ the closure of $ C^{\infty}_0(\Omega)$ in $W^{s, p(x, y)}(\Omega)$ with respect to $\|.\|_{W^{s, p(x, y)}(\Omega)}$, then  $W^{s, p(x, y)}_0(\Omega)$ is a Banach space under the Gagliardo semi-norm $[u]_{s,p(x,y)}$ given by \eqref{gagliardo}.\\    
An important role in working with the variable exponent fractional Sobolev space  is played by the
$p(\cdot,\cdot)$-modular defined in this manner:
 $$  \rho_{s,p(x,y)}^0(w)= \iint_{\Omega \times \Omega} \frac{|w(x)-w(y)|^{p(x, y)}}{|x-y|^{sp(x, y)+N}}dx dy, \quad \forall w \in Y_0.$$ 
 Similar to Proposition \ref{pro1}, we have the following results.
 \begin{proposition}[{\cite[Lemma 2.2]{azroul1}}]\label{prop2}
  Let $p: \overline{\Omega} \times \overline{\Omega} \longrightarrow(1,+\infty)$ be a continuous variable exponent that checks conditions \eqref{p1} and \eqref{p2}. Then, for any $w \in Y_0 $ and $\{w_k\}\subset Y_0$, we have
   \renewcommand{\theenumi}{\roman{enumi}}%
           \begin{enumerate}
 \item $1 \leqslant[w]_{s, p(x, y)} \Rightarrow[w]_{s, p(x, y)}^{p^{-}} \leqslant \rho_{s, p(x, y)}^{o}(w) \leqslant[w]_{s, p(x, y)}^{p^{+}}$,\\
 \item $[w]_{s, p(x, y)} \leqslant 1 \Rightarrow[w]_{s, p(x, y)}^{p^{+}} \leqslant \rho_{s, p(x, y)}^{o}(w) \leqslant[w]_{s, p(x, y)}^{p^{-}}$,\\
 \item $\lim _{k \rightarrow+\infty}\left[w_{k}-w\right]_{s, p(x, y)}=0 \Longleftrightarrow \lim _{k \rightarrow+\infty} \rho_{s, p(x, y)}^{o}\left(w_{k}-w\right)=0$.
 \end{enumerate}
 \end{proposition}
In order to correctly encode the nonlocal Dirichlet boundary datum on $\mathbb{R}^N\setminus\Omega$ in the variational formulation, the author, in \cite{azroul2},  adjusted the definition of $W^{s, p(x, y)}(\Omega)$ by extending the integral from $\Omega\times\Omega$ to  $Q:=\mathbb{R}^{2N}\setminus (\Omega^c\times \Omega^c)$ with  $\Omega^c=\mathbb{R}^N\setminus\Omega$ (see also \cite{chen, servadei} for the constant exponent case), and they gave another form of variable exponent fractional Sobolev space
$${\small Z=W^{s, p(x, y)}(Q)=~ \left\{
   \begin{array}{clcclc}
   u:\mathbb{R}^N\longrightarrow\mathbb{R} ~\text{measurable, such that}~~u_{|\Omega}\in L^{\bar{p}(x)}(\Omega) ~\text{with}~\\   \int_{Q }\tfrac{|u(x)-u(y)|^{p(x,y)}}{\lambda^{p(x,y)}|x-y|^{sp(x,y)+N}}~dxdy <+\infty,~~ \text{for some}~~\lambda>0\\
  
   \end{array}
   \right\}},
   $$ 
  where  $p:\overline{Q}\longrightarrow(1,+\infty)$ satisfies \eqref{p1} and \eqref{p2} on $\overline{Q}$.\\
   The space $W^{s, p(x, y)}(Q)$ is endowed with the following norm
   $$\|u\|_{Z}=\|u\|_{L^{\bar{p}(x)}(\Omega)}+[u]_{{s,p(x,y)}(Q)},$$ 
  where is a Gagliardo seminorm with variable exponent, defined by $$[u]_{Z}:=[u]_{s,p(x,y)}(Q)= \inf \bigg\{\lambda>0:\int_{Q}\frac{|u(x)-u(y)|^{p(x,y)}}{\lambda^{p(x,y)}|x-y|^{sp(x,y)+N}}~dxdy \leqslant1 \bigg\}.$$ 
  Similar to the space $(Y, \|.\|_{Y})$, we have that  $(Z, \|.\|_{Z})$ is a separable reflexive Banach space (see \cite[Lemma 2.3]{azroul2}). Moreover, from \cite[Theorem 2.2]{azroul2} the space $Z$ is continuously embedded in $Y$.\\
Next, let $Z_0$ denotes the following linear subspace of $Z$
  $$Z_0=\left\lbrace u\in Z: u=0 \text{ a.e. in } \mathbb{R}^N\setminus\Omega \right\rbrace,$$
  with the norm 
  $$\|u\|_{Z_0}=[u]_{Z}= \inf \bigg\{\lambda>0:\int_{Q}\frac{|u(x)-u(y)|^{p(x,y)}}{\lambda^{p(x,y)}|x-y|^{sp(x,y)+N}}~dxdy \leqslant1 \bigg\}.$$ 
 \section{Anisotropic fractional Sobolev space with variable exponent}\label{S3}
The aim of this section is to discuss some fundamental properties of the new variable exponent anisotropic fractional Sobolev space such as completeness, separability, and reflexivity. Besides, 
we give a detailed proof of the  continuous and compact embedding result. Moreover, we will adjust the definition of these spaces in order to encode correctly the Dirichlet boundary datum in nonlocal problems.
\subsection{Some qualitative properties}
\begin{theorem}\label{theorem1} The space $\left( W,\|.\|_{s,\vec{p}(x,y)} \right) $  defined by \eqref{norm} 
 is a Banach space.

\end{theorem}
\noindent{\bf {\textit{Proof}. }}
We define the operator $\mathcal{T}$ by 
$$
\begin{aligned}\mathcal{T}:W^{s, p_M(x),\vec{p}(x,y)}(\Omega) &\longrightarrow \Pi=L^{p_{M}(x)}(\Omega) \times \prod_{i=1}^{N} L^{p_{i}(x, y)}(\Omega \times \Omega)\\
w \quad & \longmapsto\left(w(x), \tfrac{w(x)-w(y)}{|x-y|^{\frac{N}{p_{1}(x, y)}+s}},\cdots, \tfrac{w(x)-w(y)}{|x-y|^{\frac{N}{p_{N}(x, y)}+s}}\right).
\end{aligned} $$
For any $w \in W$, we have 
$$ \begin{aligned}
\|\mathcal{T}(w)\|_{\Pi}  &=  \|w\|_{\small  L^{p_{M}(x)}(\Omega)} +  \sum_{i=1}^{N} \|w\|_{\small L^{p_{i}(x, y)}(\Omega \times \Omega)   }\\
&= \|w\|_{s,\vec{p}(x,y)}.  \end{aligned}$$
Then $\mathcal{T}$ is an isometry from $W$ to the Cartesian product  $\Pi$.
Since  $\Pi$ is a Banach space, therefore $W$ is also a Banach space.
\begin{proposition} Let $\Omega$ be
 a bounded open domain  of $\mathbb{R}^N$, we suppose that conditions
\eqref{p1} and \eqref{p2} are satisfied. Then $\left( W,\|.\|_{s,\vec{p}(x,y)} \right) $ is a separable and uniformly convex space $($consequently it is a reflexive space$)$. 
\end{proposition}
\noindent{\bf {\textit{Proof}. }}
We use the same isometric map introduced in the proof of Theorem \ref{theorem1}. Because $W$ is a Banach space, we get that $\mathcal{T}\left( W\right) $ is a closed subset of $\Pi$. Consequently, by \cite[Proposition 3.20 and Proposition 3.25]{brezis}, it follows that the space $\left( W,\|.\|_{s,\vec{p}(x,y)} \right) $ is separable.\\
On the other hand, since $p_i^->1$, for $i=1,\cdots,N$ from \cite{Fan} the space $(L^{p_i(x)},\|.\|_{p_i(x)})$ is  uniformly convex. Moreover, by a theorem in \cite[p. 184]{Musielak}, we infer that $\Pi$ is uniformly convex space.  As $\mathcal{T}\left( W\right) $ is a closed subset of $\Pi$  and since any linear subspace of a uniformly convex linear normed space is also uniformly convex, hence $\left( W,\|.\|_{s,\vec{p}(x,y)} \right) $ is uniformly convex and consequently it is a reflexive space.
 \\
 
Now, let us come back to the definition of the space $W=W^{s, p_M(x),\vec{p}(x,y)}(\Omega)$. Since  $\Omega\times\Omega$ is strictly contained in $Q$ and for the same reason mentioned in Section \ref{S2} on the definition of the space $Y= W^{s, p(x,y)}(\Omega)$, we will adjust    the definition of the space $W$
as follows 
 $$
 \begin{aligned}
  X&=W^{s, p_M(x),\vec{p}(x,y)}(Q)\\
  &=\left\{
   \begin{array}{clcclc}
   w:\mathbb{R}^N\longrightarrow\mathbb{R} ~\text{measurable, such that}~~w_{|\Omega}\in  L^{p_M(x)}(\Omega) ~\text{with}~\\   \int_{Q}\tfrac{\left| w(x)-w(y)\right|^{p_i(x,y)} }{\nu_i^{p_i(x,y)}\left| x-y \right| ^{sp_i(x,y)+N} } dx dy <+\infty,~~  \text{for  some}~~\nu_i>0 \text{ and for all }i=1,...,N
   \end{array}
   \right\},
\end{aligned} $$
where $\vec{p}:\overline{Q}^N \longrightarrow \mathbb{R}^N$ is a continuous vectoriel function with $p_i$ satisfying conditions \eqref{p1} and \eqref{p2} on $\overline{Q}$.
\begin{remark}\text{}
\begin{enumerate}
\item[$(i)$-] Simply swap $\Omega \times \Omega$ for $Q$ to define the norm on $X$ in the same manner on $W$, and we denote this norm by $\|.\|_{X} $.
\item[$(ii)$-]  Similar to the space $(W, \|.\|_{W})$, we have that  $\left( X,\|.\|_{X} \right) $ is a Banach, reflexive, and separable space.
\end{enumerate}
\end{remark} 
We define  on $X$ the following mapping 
$$ \rho_{s,\vec{p}(x, y)}(w)=\int_{\Omega}|w(x)|^{p_M(x)} d x + \sum_{i=1}^{N}\int_{Q} \frac{|w(x)-w(y)|^{p_{i}(x, y)}}{|x-y|^{sp_{i}(x, y)+N}}dx dy, $$
for all $w\in X$, called $\vec{p}(x, y)$-modular. It is obvious to check that $ \rho_{\vec{p}(., .)}$ is a convex modular on $X$. Moreover, the modular norm on $X$ defined by $$ \|w\|_{ \rho_{s,\vec{p}}}= \inf  \left\lbrace \nu  >0, \quad  \rho_{s,\vec{p}\left( x, y\right) }\left( \frac{w}{\nu}\right) \leqslant 1 \right\rbrace    $$
is equivalent to the norm $\|.\|_{X}$.\\
Denote by $X_0$ the closed subspace of $X$ introduced as follows
$$
X_0=W_{0}^{s, p_M(x),\vec{p}(x,y)}(Q)=\left\{w \in X: w(x)=0 \text { a.e. } x \in \mathbb{R}^{N} \backslash \Omega\right\}.
$$ 
 This subset of $X$ also qualifies as a Banach space which is separable and reflexive.\\
 Besides, for any $w \in  X_0$, we define the convex modular as follows
  $$ \rho_{s,\vec{p}(x, y)}^0(w)= \sum_{i=1}^{N}\int_{Q} \frac{|w(x)-w(y)|^{p_{i}(x, y)}}{|x-y|^{sp_{i}(x, y)+N}}dx dy. $$
  The modular on $X_0$ verifies the following properties.
  \begin{lemma} \label{l2.1}
  Let $\vec{p}:\overline{Q}^N \longrightarrow \left( 1,+\infty\right) $ be a vector function satisfying conditions \eqref{p1} and \eqref{p2}. Then, for any $w \in X_0$, we have 
  \begin{itemize}
  \item[$(i)$-]  $ 1 \leqslant[w]_{s, \vec{p}(x, y)} \Rightarrow \rho_{s, \vec{p}(x,y)}^0(w) \leqslant[w]_{s, \vec{p}(x, y)}^{P^{+}_+},$ \\
  \item[$(ii)$-] $ [w]_{s, \vec{p}(x, y)} \leqslant 1 \Rightarrow\rho_{ s,\vec{p}(x,y)}^0(w)  \leqslant [w]_{s, \vec{p}(x, y)}.$
  \end{itemize}
  \end{lemma}
\noindent{\bf{\textit{Proof}.}}
   $(i)$- Let $w \in X_0$ such that $[w]_{s, \vec{p}(x, y)} \geqslant 1$, then
$\sum_{i=1}^{N}[w]_{s, p_i(x, y)} \geqslant 1 $. Hence, 
$$ \exists J\subset I =\left\lbrace 1,...,N\right\rbrace \text{ such that } \forall j \in J \quad [w]_{s, p_j(x, y)} \leqslant 1 \text{ and }  \forall i \in I\setminus J \quad [w]_{s, p_i(x, y)} \geqslant 1,$$ 
which implies by using Proposition \ref{prop2}
$$
\begin{aligned}
\rho_{ s,\vec{p}(x,y)}^0(w) &= \sum_{j\in J} \rho_{s, p_j(x,y)}^0(w) + \sum_{i\in I\setminus J } \rho_{s, p_i(x,y)}^0(w)\\
& \leqslant \sum_{j\in J} [w]_{s, p_j(x, y)}^{p_j^-} + \sum_{i\in I\setminus J } [w]_{s, p_i(x, y)}^{p_i^+}\\
& \leqslant \sum_{j\in J} [w]_{s, p_j(x, y)} + \left( \sum_{i\in I\setminus J } [w]_{s, p_i(x, y)}\right)^{P^+_+}\\
& \leqslant [w]_{s, \vec{p}(x, y)}^{P^{+}_+}
\end{aligned}
.$$
$(ii)$- Let $w \in X_0$ such that $[w]_{s, \vec{p}(x, y)} \leqslant 1$, thus $\forall i=1,...,N \quad [w]_{s, p_i(x, y)} \leqslant 1$,
consequently, $  \rho_{s,\vec{p}(x, y)}^0(w) \leqslant \sum_{i=1}^{N} [w]_{s, p_i(x, y)}^{p_i^-}  \leqslant  [w]_{s, \vec{p}(x, y)}$.
 
\begin{remark}
The results of Lemma $\ref{l2.1}$ can be extend to the modular $\rho_{s,\vec{p}(., .)}$ defined on $X$.
\label{rem2.1}
\end{remark}
\subsection{Embedding Results} In the present subsection, we focus on the continuous and compact embedding of our new functional framework into Lebesgue spaces. We start by giving a detailed proof of Theorem \ref{th1} which is essential in the study of nonlocal anisotropic problems. \\
\textbf{Proof of Theorem  \ref{th1}} We have $\vec{p}, p_M$, and $q$ are continuous variable exponents, then  for $i=1,\ldots, N$, we can find a family of positive constants $\varepsilon^1_i$, and $\varepsilon^2_i$ such that    
\begin{equation}  
\left\{\begin{array}{rlrl}
\frac{N p_i(x, x)}{N-s p_i(x, x)}-q(x)& \geqslant \frac{\varepsilon^1_i}{2}, \\
\quad p_M(x)-p_i(x, x)& \geqslant \frac{\varepsilon^2_i}{2},
\end{array} \right.
 \end{equation}
 for all $x \in \overline{\Omega}$.\\ On the other hand, by using the finite covering theorem, we can find a positive constant $\varepsilon=\varepsilon \left(\vec{p}, p_M,q,\varepsilon^1_i,\varepsilon^2_i,t  \right),$ and a finite family of open sets $B_k \subset \Omega$ such that $ \Omega= \bigcup _{k=1}^N B_k $ with $ diam(B_k)< \varepsilon$, and $t\in (0,s)$. Hence, for $i=1,\ldots, N$, we have
 \begin{equation}  
 \left\{\begin{array}{rlrl}
 \frac{N p_i(y, z)}{N-t p_i(y, z)}-q(x)& \geqslant \frac{\varepsilon^1_i}{2}, \\
 \quad p_M(x)-p_i(y, z)& \geqslant \frac{\varepsilon^2_i}{2},
 \end{array} \right.
 \label{3}
  \end{equation} 
  for all $(x,y,z) \in B_k^3$.\\
  Set $p_k^i=\inf _{(y, z) \in B_k \times B_k } \left\lbrace p_i(y, z)- \mu \right\rbrace $, based on \eqref{3}, we can choose $\mu >0$ small enough such that $ p_i^{-}-1> \mu>0$, for $i=1,\ldots, N$, where $$
 \frac{N p_k^i  }{N-t p_k^i }-q(x) \geqslant \frac{\varepsilon^1_i}{2}  \quad  \quad \text{and for all }x \in B_k.
  $$
  We get, for $i=1,\ldots, N$  and   all $  x \in B_k$
  \begin{enumerate}
  \item[$\bullet$] If $p_k^{i,*}=  \frac{N p_k^i  }{N-t p_k^i }$, then $p_k^{i,*} \geqslant \frac{\varepsilon^1_i}{2} +q(x), $\\
  \item[$\bullet$]   $ p_{M}(x)\geqslant \frac{\varepsilon^2_i}{2} + p_k^i.  $
  \end{enumerate}
  Due to the embedding theorem in the case of constant exponent (\cite[Theorem 6.7]{dinezza}), there exists a constant $C=C(N,t,p_k^i, B_k,\varepsilon)$ such that 
  $$
  \|w\|_{L^{p_k^{i,*}}(B_k)} \leq C\|w\|_{W^{t,p_k^i }(B_k)}.
  $$
 To prove inequality \eqref{ine1}, we need to show the following results
 \begin{itemize}
 \item[$(i)$-] $\exists C>0 $; $\quad \sum_{k=1}^N\|w\|_{L^{p_k^{i,*}}\left(B_k\right)} \geq C\|w\|_{L^{q(x)}(\Omega)}.$\\
 \item[$(ii)$-]  $\exists C>0 $; $\quad  C\|w\|_{L^{p_{M}(x)}(\Omega)} \geq \sum_{k=1}^N\|w\|_{L^{p_k^i}\left(B_k\right)}.  $\\
 \item[$(iii)$-]  $\exists C>0 $; $\quad C\left[ w\right]_{s, p_i(x, y),\Omega}  \geq \sum_{k=1}^N[w]_{t, p_k^i,\left(B_k\right)}$.
 \end{itemize}
 Note that the constant $C$ may differ from line to line.\\
 $(i)$- Indeed, we have $|w(x)|=\sum_{k=1}^N |w(x)| \chi_{B_k}(x) $, then  \begin{equation}
 \|w\|_{L^{q(x)}(\Omega)} \leq \sum_{k=1}^N\|w\|_{L^{q(x)}\left(B_k\right)}.
 \label{5}
 \end{equation}
 Since $p_k^{i,*}> q(x),$ for all $x \in B_k$. Hence, we can take $\alpha_k^i \in C^+(B_k)$ such that $\frac{1}{q(x)}=\frac{1}{p_k^{i,*}}+\frac{1}{\alpha_k^i(x)}$, for all $x\in B_k$, from Corollary \ref{H2.1}, we get 
\begin{equation}\begin{aligned}
 \|w\|_{L^{q(x)}\left(B_k\right)} & \leq C\|w\|_{L^{p_k^{i,*}(x)}\left(B_k\right)}\|\chi_{B_k}\|_{L^{\alpha_k^i(x)}\left(B_k\right)} \\
 & \leq C\|w\|_{L^{p_k^{i,*}(x)}\left(B_k\right)} .
 \end{aligned}    \label{6}
  \end{equation}
  Combining \eqref{5}, and \eqref{6} we obtain the first inequality.\\
 $(ii)$- From the fact that $p_{M}(x)>p_k^{i},$ for all $x \in B_k$. then there exists $\beta_k^i \in C^+(B_k)$ such that $\frac{1}{p_k^{i}}=\frac{1}{p_{M}(x)}+\frac{1}{\beta_k^i(x)}$. Following, we proceed as before, we get the desired inequality.\\
 $(iii)$-  We consider $ \mathcal{D}(x,y)=\frac{|w(x)-w(y)|}{|x-y|^s}  $ and for 	all $  i=1,\ldots, N $,  we set $d\eta_k(x,y)= \frac{dx dy }{|x-y|^{N+(t-s)p_k^i}}$  , where $\eta_k$ represents a measure in $B_k \times B_k$ with $\eta_k (B_k \times B_k) < \infty$. We have 
 $$ \begin{aligned}
   [w]_{t, p_k^i, B_k}&=\left(\int_{B_k} \int_{B_k} \frac{|w(x)-w(y)|^{p_k^i}}{|x-y|^{N+t p_k^i}} \mathrm{~d} x \mathrm{~d} y\right)^{\frac{1}{p_k^i}}\\
 & =\left(\int_{B_k} \int_{B_k}\left(\frac{|w(x)-w(y)|}{|x-y|^s}\right)^{p_k^i} \frac{\mathrm{d} x \mathrm{~d} y}{|x-y|^{N+(t-s) p_k^i}}\right)^{\frac{1}{p_k^i}} \\
 & =\|\mathcal{D}\|_{L^{p_k^i}\left(\eta_k, B_k \times B_k\right)}.
 \end{aligned}   $$ 
Using H\"{o}lder's inequality, we get $$
\|\mathcal{D}\|_{L^{p_k^i}\left(\eta^k, B_k \times B_k\right)}\leqslant C \|\mathcal{D}\|_{L^{p_i(x,y)}\left(\eta_k, B_k \times B_k\right)}.
$$
For  $   i=1,\ldots, N  $, take $\lambda_i>0$, such that $$ \int_{B_k} \int_{B_k} \frac{|w(x)-w(y)|^{p_i(x,y)}}{\lambda_i^{p_i(x,y)}|x-y|^{N+t p_i(x,y)}} \mathrm{~d} x \mathrm{~d} y <1,$$ 
and choose $m:=\max \left\{1, \sup _{(x, y) \in \Omega \times \Omega}|x-y|^{s-t}\right\}$, and $\tilde{\lambda}_i=m \lambda_i$, we have
$$
\begin{aligned}
  \int_{B_k} \int_{B_k}\left(\tfrac{|w(x)-w(y)|}{\tilde{\lambda_i}|x-y|^s}\right)^{p_i(x, y)} \mathrm{d} \eta_k(x, y) 
&=\int_{B_k} \int_{B_k}\left(\tfrac{|w(x)-w(y)|}{\tilde{\lambda_i}|x-y|^s}\right)^{p_i(x, y)} \tfrac{\mathrm{d}x \mathrm{d}y }{|x-y|^{N+(t-s)p_i^k}}  \\
& =\int_{B_k} \int_{B_k} \tfrac{|w(x)-w(y)|^{p_i(x, y)}}{(\lambda_i )^{p_i(x, y)}|x-y|^{N+s p_i(x, y)}} \tfrac{|x-y|^{(s-t) p_i^k}}{m^{p_i(x,y)}} \mathrm{d} x \mathrm{d} y \\
& \leq \int_{B_k} \int_{B_k} \frac{|w(x)-w(y)|^{p_i(x, y)}}{(\lambda_i )^{p_i(x, y)}|x-y|^{N+s p_i(x, y)}}\mathrm{d} x \mathrm{~d} y <1 .
\end{aligned}
$$
Then, we get  $$ \|\mathcal{D}\|_{L^{p_i(x,y)}\left(\eta_k, B_k \times B_k\right)} <  \tilde{\lambda}_i.$$
Consequently,  $$ \|\mathcal{D}\|_{L^{p_i(x,y)}\left(\eta_k, B_k \times B_k\right)} \leq m \left[ w\right]_{s, p_k^i,\Omega}~~\Longrightarrow ~~ \left[ w\right]_{t, p_k^i,B_k} \leq C\left[ w\right]_{s, p_k^i,\Omega}.$$  
Now, from  the  inequalities $(i)$-$(iii)$, we infer 
\begin{equation}
\begin{aligned}
& \|w\|_{L^{q(x)}}(\Omega) \leq C \sum_{k=1}^{N}\|u\|_{L^{p_k^{i, *}}{ }_{\left(B_k\right)}}  \\
& \leq C\sum_{k=1}^{N}\left( \|w\|_{L^{p_k^i}\left(B_k\right)}+[w]^{t, p_k^i}\left(B_k\right)\right) \\
& \leq C\left(\|w\|_{L^{p_M(x)} (\Omega)}+[w]^{s, p_i(x, y)}(\Omega)\right)  \\
& \leq C\left(\|w\|_{L^{p_M(x)} (\Omega)}+\sum_{k=1}^{N}[w]^{s, p_i(x, y)}(\Omega)\right) \\
&= C  \|w\|_{s,\vec{p}(x,y)}.
\end{aligned}
\end{equation} 
Next, we investigate the compactness of the embedding between $W$ and $L^{q(x)} (\Omega)$.
Let $\left\lbrace w_n\right\rbrace $ be a bounded sequence on $W$, that is, there exists $M>0$ such that $\|w_n\|_{s,\vec{p}(x,y)} \leq M$.
It follows that $ [w_n]_{t, p_k^i \left(B_k\right)} \leq M   $, and by the compact embedding in the constant case, we can find a subsequence of $\left\lbrace w_n\right\rbrace $ such that  for $k=1,...,N$, $\{w_n^k\}$ converges strongly in $L^{p_k^{i,*}-\frac{\varepsilon}{3}}(B_k)$.\\
Since $|w(x)|=\sum_{k=1}^N |w(x)| \chi_{B_k}(x) $, therefore
$$  
\|w_n^k - w \|_{L^{q(x)}(\Omega)} \leqslant   C\sum_{k=1}^N\|w_n^k -w \|_{L^{p_k^{i,*}}\left(B_k\right)}.
$$
As a consequence, we obtain the strong convergence in $L^{q(x)}(\Omega)$. Using the same procedures as the previous proof we get the desired results.
\begin{theorem}
Consider a smooth bounded domain $\Omega$ in $\mathbb{R}^N$ $(N\geqslant 2)$, and let $\vec{p}:\overline{Q}^N \longrightarrow \left( 1,+\infty\right) $ be a vector function satisfying conditions \eqref{p1} and \eqref{p2}. Moreover, assume that for $i=1...N$,  $sp_i(x,y)<N$ for all $ (x,y)\in \overline{Q}$. Then, we have :
\begin{itemize}
\item[$(i)$-] $W$ is continuously embedded in $X$.
\item[$(ii)$-] The results of Theorem $\ref{th1}$ hold true for the space $X$ and $X_0$.
\end{itemize}
\label{th2}
\end{theorem}
\noindent{\bf{\textit{Proof}.}} 
  $(i)$- Let $ w \in X$, since $\Omega \times \Omega \subsetneq Q$, then for $  i=1,...,N$, and for some $\nu_i >0$ 
$$
\int_{\Omega \times \Omega} \frac{|w(x)-w(y)|^{p_i(x, y)}}{\nu_i^{p_i(x, y)}|x-y|^{s p_i(x, y)+N}} \mathrm{~d} x \mathrm{~d} y \leqslant \int_Q \frac{|w(x)-w(y)|^{p_i(x, y)}}{\nu^{p_i(x, y)}|x-y|^{s p_i(x, y)+N}} \mathrm{~d} x \mathrm{~d} y .$$
As a consequence, for $  i=1,...,N$, we get  
$${\footnotesize \tiny \left\{
\nu_i>0: \int_{Q} \tfrac{|w(x)-w(y)|^{p_i(x, y)}}{\nu^{p_i(x, y)}|x-y|^{s p_i(x, y)+N}} \mathrm{~d} x \mathrm{~d} y \leqslant 1  \right\}  \subset 
\left\{\nu_i>0: \int_{\Omega \times \Omega} \tfrac{|w(x)-w(y)|^{p_i(x, y)}}{\nu^{p_i(x, y)}|x-y|^{s p_i(x, y)+N}} \mathrm{~d} x \mathrm{~d} y \leqslant 1\right\} }.$$ 
Which implies that 
$$
\begin{aligned}
\left[ w\right]_{s,\vec{p}(x,y),(\Omega)}=& \sum_{i=1}^{N} \inf \left\{\nu_{i}>0: \iint_{\Omega \times \Omega} \frac{|w(x)-w(y)|^{p_{i}(x, y)}}{\nu_{i}^{p_{i}(x, y)}|x-y|^{sp_{i}(x, y)+N}}<1\right\} \\
 \leqslant&  \sum_{i=1}^{N} \inf \left\{\nu_{i}>0: \int_{Q} \frac{|w(x)-w(y)|^{p_{i}(x, y)}}{\nu_{i}^{p_{i}(x, y)}|x-y|^{sp_{i}(x, y)+N}}<1\right\} = \left[ w\right]_{s,\vec{p}(x,y),(Q)}.
\end{aligned}
$$
Therefore, $  \|w\|_{W} \leqslant  \|w\|_{X} $, providing to prove the desired result.\\
$(ii)$-  Since $W$ is continuously embedded in $X$, it is evident that the results of Theorem \ref{th1} can be extend to the space $X$ and $X_0$.
 
\section{Basic properties of the fractional $\vec{p}(x,.)$-Laplacian functional}\label{S4} Let us introduce the following mappings which are crucial to deal with the fractional $\vec{p}(x,.)$-Laplacian equations.
For any $ w \in X$, we consider
$$
L_{\vec{p}(x,y)} \left( w\right) =\sum_{i=1}^N \int_{Q} \frac{1}{p_i(x, y)} \frac{|w(x)-w(y)|^{p_i(x, y)}}{|x-y|^{s p_i(x, y)+N}} d x d y + \int_{\Omega} \frac{1}{p_M(x)} |w(x)|^{p_M(x)} dx.
$$
In addition, for all $ v, w \in X$, we define 
$$
\begin{aligned}
\left\langle \mathcal{L}_{\vec{p}(x,y)} \left( v\right), w\right\rangle= \sum_{i=1}^N & \int_{Q} \frac{\left|v(x)-v(y)\right|^{p_{i}(x,y)-2}\left( v(x)-v(y)\right) \left( w(x)-w(y)\right) }{|x-y|^{s p_i(x, y)+N}} d x d y \\
+& \int_{\Omega}|v(x)|^{p_M(x)-2} v(x) w(x) d x.
\end{aligned}
$$
As well as, for all $ v,w  \in X_0$, we  consider $I_{\vec{p}(.,.)}$, and $\mathcal{I}_{\vec{p}(.,.)}$ in this way 
$$\begin{aligned}
I_{\vec{p}(x,y)} \left( w\right) &=\sum_{i=1}^N \int_{Q} \frac{1}{p_i(x, y)} \frac{|w(x)-w(y)|^{p_i(x, y)}}{|x-y|^{s p_i(x, y)+N}} d x d y,\\
\left\langle \mathcal{I}_{\vec{p}(x,y)} \left( v\right), w\right\rangle&= \sum_{i=1}^N \int_{Q} \frac{\left| v(x)-v(y)\right|^{p_{i}(x,y)-2}\left( v(x)-v(y)\right) \left( w(x)-w(y)\right) }{|x-y|^{s p_i(x, y)+N}} d x d y .
\end{aligned}
$$

\begin{lemma}
Suppose that hypotheses \eqref{p1} and \eqref{p2} are satisfied, then the following assertions hold true:
\begin{itemize}
\item[$(i)$-] $L_{\vec{p}(x,y)}$ and $I_{\vec{p}(x,.)}$ are well defined on $X$ and $X_0$ respectively.\\
\item[$(ii)$-]  $L_{\vec{p}(x,y)} \in \mathcal{C}^1\left(X, \mathbb{R}\right)$ and for all $v, w\in X$, we have $$\left\langle L_{\vec{p}(x,y)}^{\prime} (v), w\right\rangle = \left\langle \mathcal{L}_{\vec{p}(x,y)} \left( v\right) , w\right\rangle.$$
\item[$(iii)$-]  $I_{\vec{p}(x,y)} \in \mathcal{C}^1\left(X_0, \mathbb{R}\right)$ and for all $v, w\in X_0$, we have $$\left\langle I_{\vec{p}(x,y)}^{\prime} (v), w\right\rangle = \left\langle \mathcal{I}_{\vec{p}(x,y)} \left( v\right) , w\right\rangle .$$
\end{itemize}
\label{lemma 2.3}
\end{lemma}
\noindent{\bf{\textit{Proof}.}} $(i)$- Let $w \in X$, then by using Remark \ref{rem2.1} and Lemma \ref{l2.1}, we get
$$
\begin{aligned}
L_{\vec{p}(x,y)} \left( w\right) &=\sum_{i=1}^N \int_{Q} \frac{1}{p_i(x, y)} \frac{|w(x)-w(y)|^{p_i(x, y)}}{|x-y|^{s p_i(x, y)+N}} d x d y + \int_{\Omega} \frac{1}{p_M(x)} |u(x)|^{p_M(x)} dx\\
& \leqslant \frac{1}{P_{-}^-} \sum_{i=1}^N \int_{Q}  \frac{|w(x)-w(y)|^{p_i(x, y)}}{|x-y|^{s p_i(x, y)+N}} d x d y +  \frac{1}{p_M^-} \int_{\Omega} |w(x)|^{p_M(x)} dx \\
& \leqslant \max \left\lbrace   \tfrac{1}{P_{-}^-}, \tfrac{1}{p_M^-}   \right\rbrace \rho_{s,\vec{p}(x, y)}(w) \\
& \leqslant \max \left\lbrace   \tfrac{1}{P_{-}^-}, \tfrac{1}{p_M^-}   \right\rbrace \|w\|_{X}^{P_+^+} <+\infty.
\end{aligned}
$$
Consequently, $L_{\vec{p}(x,y)}$ is well defined on $X$. By similar argument, $I_{\vec{p}(x,y)}$ is also well defined on $X_0 $. \\
$(ii)$ and $(iii)$- Standard reasoning as in \cite{azroul1} yields the second and third results.
 
\begin{lemma}
Under conditions \eqref{p1} and \eqref{p2}, we have:
\begin{itemize}
\item[$(i)$-] The operator $\mathcal{I}_{\vec{p}(.,.)}$ is  strictly monotone and bounded.\\
\item[$(i)$-] The operator $\mathcal{I}_{\vec{p}(.,.)}$ is of type $(S_+)$.
\end{itemize}
\end{lemma}
\noindent{\bf{\textit{Proof}.}}
$(i)$-  For any $v, w \in X_0$ with $v \neq w$, we have by using Cauchy's inequality: 
{\small $$
\begin{aligned}
\left\langle \mathcal{I}_{\vec{p}(x,y)}(v)- \mathcal{I}_{\vec{p}(x,y)}(w),v-w \right\rangle & = \left\langle \mathcal{I}_{\vec{p}(x,y)}(v),v\right\rangle - \left\langle \mathcal{I}_{\vec{p}(x,y)}(v),w\right\rangle - \left\langle \mathcal{I}_{\vec{p}(x,y)}(w),v\right\rangle + \left\langle  \mathcal{I}_{\vec{p}(x,y)}(w),w \right\rangle \\
& \geqslant \frac{1}{2} \sum_{i=1}^N \int_{Q} \tfrac{|v(x)-v(y)|^{p_i(x, y)-2} \left( \left| v(x)-v(y)\right| ^2-\left| w(x)-w(y)\right|^2\right) }{|x-y|^{s p_i(x, y)+N}} d x d y \\ 
& - \frac{1}{2} \sum_{i=1}^N \int_{Q} \tfrac{|w(x)-w(y)|^{p_i(x, y)-2} \left( \left| v(x)-v(y)\right| ^2-\left| w(x)-w(y)\right|^2\right) }{|x-y|^{s p_i(x, y)+N}} d x d y.
\end{aligned}
$$ }
By changing the role of $v$ and $w$ and the Simon's inequality \cite{simon}, we get
{\small  $$
 \begin{aligned}
&\left\langle \mathcal{I}_{\vec{p}(x,y)}(v)- \mathcal{I}_{\vec{p}(x,y)}(w),v-w \right\rangle \\
&\geqslant  \tfrac{1}{2} \sum_{i=1}^N \int_{Q} \tfrac{\left(  \left| v(x)-v(y)\right| ^{p_i(x, y)-2}- \left| w(x)-w(y)\right| ^{p_i(x, y)-2} \right)  \left( \left| v(x)-v(y)\right| ^2-\left| w(x)-w(y)\right|^2\right) }{|x-y|^{s p_i(x, y)+N}} d x d y   > 0 .
 \end{aligned}
$$}
Which implies that $ \mathcal{I}_{\vec{p}(.,.)}$ is strictly monotone. It is self-evident to demonstrate the boundedness of $ \mathcal{I}_{\vec{p}(.,.)}$.\\
$(ii)$- 
Consider a sequence $\{w_n\}$ such that $w_n \rightharpoonup w$ in $X_0$ and \\$\limsup _{n \rightarrow+\infty}\left\langle \mathcal{I}_{\vec{p}(x,y)}\left(w_n\right)-\mathcal{I}_{\vec{p}(x,y)}(w), w_n-w\right\rangle \leq 0$, from the strict monotony of $\mathcal{I}_{\vec{p}(.,.)}$, we obtain  that 
\begin{equation}
\limsup _{n \rightarrow+\infty}\left\langle \mathcal{I}_{\vec{p}(x,y)}\left(w_n\right)-\mathcal{I}_{\vec{p}(x,y)}(w), w_n-w\right\rangle = 0.
\end{equation}
In the view of Theorem \ref{theorem1}, $w_n(x) \longrightarrow w(x)$ a.e in $\Omega$.
By using Fatou's Lemma, we have:
{\begin{equation}
\small  \liminf _{n \rightarrow+\infty} \sum_{i=1}^N \int_{Q}  \frac{|w_n(x)-w_n(y)|^{p_i(x, y)}}{|x-y|^{s p_i(x, y)+N}} d x d y \geqslant \sum_{i=1}^N \int_{Q} \frac{|w(x)-w(y)|^{p_i(x, y)}}{|x-y|^{s p_i(x, y)+N}} d x d y.
\label{fatou}
\end{equation}}
Young's inequality implies the existence of two constants $C_1$ and $C_2$ such that
{\small \begin{equation}
\begin{aligned}
\left\langle \mathcal{I}_{\vec{p}(x,y)}\left( w_n\right), w_n-w\right\rangle &= \left\langle \mathcal{I}_{\vec{p}(x,y)}\left(w_n\right), w_n\right\rangle - \left\langle \mathcal{I}_{\vec{p}(x,y)}\left(w_n\right), w\right\rangle \\
& \geqslant \sum_{i=1}^N \int_{Q}  \frac{|w_n(x)-w_n(y)|^{p_i(x, y)}}{|x-y|^{s p_i(x, y)+N}} d x d y \\
&- \sum_{i=1}^N \int_{Q}  \left( \tfrac{p_i(x,y)-1}{p_i(x,y)} \tfrac{|w_n(x)-w_n(y)|^{p_i(x, y)}}{|x-y|^{s p_i(x, y)+N}} + \tfrac{1}{p_i(x,y)} \tfrac{w(x)-w(y)|^{p_i(x, y)}}{|x-y|^{s p_i(x, y)+N}}  \right)  d x d y \\
& \geqslant C_1 \sum_{i=1}^N \int_{Q}  \tfrac{|w_n(x)-w_n(y)|^{p_i(x, y)}}{|x-y|^{s p_i(x, y)+N}} d x d y - C_2 \sum_{i=1}^N \iint_{Q}  \tfrac{|w(x)-w(y)|^{p_i(x, y)}}{|x-y|^{s p_i(x, y)+N}} d x d y.
\end{aligned}
\label{young}
\end{equation} }
Combining \eqref{fatou} and \eqref{young}, we obtain 
\begin{equation}
\small  \lim _{n \rightarrow+\infty} \sum_{i=1}^N \int_{Q}  \frac{|w_n(x)-w_n(y)|^{p_i(x, y)}}{|x-y|^{s p_i(x, y)+N}} d x d y = \sum_{i=1}^N \int_{Q} \frac{|w(x)-w(y)|^{p_i(x, y)}}{|x-y|^{s p_i(x, y)+N}} d x d y.
\end{equation}
According to \cite{Natanson}, and the proof of \cite[Lemma 3.1 ]{ayazoglu}, we obtain the strong convergence of $w_n$ to $w$ in $X_0$.
  
\section{An application to anisotropic nonlocal problems}\label{S5}
As an application of this framework
we consider a class of fractional anisotropic problems with variable exponent of the following form
$$
 \left( \mathcal{\vec{P}}\right) 
 \left\{\begin{array}{ll} \left( -\Delta\right) ^s_{\vec{p}(x,.)} \left( u\right) + \left| u\right|^{p_M(x)-2}u= \lambda f(x,u) & \text { in } \Omega
, \\
\hspace{3cm} u=0 & \text { on } \mathbb{R}^{N} \backslash \Omega.
 \end{array}\right.
 \label{eq} 
$$
Where 
\begin{itemize}
\item[$\bullet$]  $\Omega $ is a bounded domain of $\mathbb{R}^N$ with Lipschitz boundary $(N \geqslant 2)$. 
\item[$\bullet$]   $\lambda >0$ is a real parameter.
\item[$\bullet$] $f:\Omega \times \mathbb{R}\longrightarrow \mathbb{R}$ is a continuous function that assumes throughout this section equal to  $\left| u\right|^{r(x)-2}u$, with 
$$
1<r^{-} \leq r(x)\leq r^{+}< P^-_- \leq p_M^{+} < +\infty, \quad \forall x \in \Omega.
$$
\end{itemize}

In  \cite{mihalescu, mihalescu2} the authors investigated the integer case of problem \hyperref[eq]{$ \left( \mathcal{\vec{P}}\right)$} and in the fractional scenario, this work \cite{azroul1} is well recommended.\\

\noindent A function $w \in X_0$ is defined as a weak solution of problem  \hyperref[eq]{$ \left( \mathcal{\vec{P}}\right)$} if :
\begin{equation}
\begin{gathered}
\sum_{i=1}^N \int_{Q} \frac{|w(x)-w(y)|^{p_i(x, y)-2}(w(x)-w(y))(\varphi(x)-\varphi(y))}{|x-y|^{s p_i(x, y)+N}} d x d y \\
+\int_{\Omega}|w(x)|^{p_M(x)-2} w(x) \varphi(x) d x-\lambda \int_{\Omega}|w(x)|^{r(x)-2} w(x) \varphi(x) d x, \quad \forall \varphi \in X_0 .
\end{gathered}
\end{equation}
For all $w\in X_0$, the energy functional related to problem \hyperref[eq]{$ \left( \mathcal{\vec{P}}\right)$}  is of the following form 
$$
\begin{aligned}
\mathcal{J}_{\lambda}(w)&= L_{\vec{p}(x,y)}(w)- \lambda \int_{\Omega} \frac{1}{r(x)}  |w(x)|^{r(x)} dx\\
&= \sum_{i=1}^N \int_{Q} \tfrac{1}{p_i(x, y)} \tfrac{|w(x)-w(y)|^{p_i(x, y)}}{|x-y|^{s p_i(x, y)+N}} d x d y +  \int_{\Omega}
 \tfrac{|w(x)|^{p_M(x)}}{p_M(x)}  dx - \lambda \int_{\Omega} \tfrac{|w(x)|^{r(x)}}{r(x)}   dx .
\end{aligned}
$$
By a standard argument as in Lemma \ref{lemma 2.3}, we have that $\mathcal{J}_{\lambda} \in C^1\left(X_0, \mathbb{R} \right)$, in addition, its derivative is given by the following formula :
\begin{equation}
\begin{aligned}
\left\langle \mathcal{J}_{\lambda}^{\prime}(w),v\right\rangle  = &\sum_{i=1}^N \int_{Q} \frac{|w(x)-w(y)|^{p_i(x, y)-2}(w(x)-w(y))(v(x)-v(y))}{|x-y|^{s p_i(x, y)+N}} d x d y \\
&+\int_{\Omega}|w(x)|^{p_M(x)-2} w(x) v(x) d x-\lambda \int_{\Omega}|w(x)|^{r(x)-2} w(x) v(x) d x,
\end{aligned}
\end{equation} 
for all $w,v \in X_0$. Hence a critical point of $\mathcal{J}_{\lambda}$ coincides with weak solution of problem \hyperref[eq]{$ \left( \mathcal{\vec{P}}\right)$}.\\
The result of the existence of solutions of problem \hyperref[eq]{$ \left( \mathcal{\vec{P}}\right)$} is formulated in the next theorem.
\begin{theorem}\label{th3.1} Let  $\Omega$ be a smooth
 bounded domain  in  $\mathbb{R^N}$ and $s \in (0,1)$. Let $\vec{p}:\overline{Q}^N \longrightarrow \left( 1,+\infty\right) $ be a vector function which supposed to verify conditions \eqref{p1}, and \eqref{p2} on $\overline{Q}$. Moreover, consider $r:\Omega \longrightarrow \left( 1, +\infty \right) $ as a continuous variable exponent such that \begin{equation}
 1< r(x) <  P^-_-  \quad \text{for all } x\in \overline{\Omega}.
 \label{condition1}
\end{equation}
Then, there exists $\lambda_* >0$ such that problem \hyperref[eq]{$ \left( \mathcal{\vec{P}}\right)$} possesses a weak solution $\forall \lambda \in (0, \lambda_*)$.
\label{existence}
\end{theorem}
In order to establish the proof of Theorem \ref{th3.1}, we will apply the mountain pass theorem combined with the Ekeland's variational principle. The crucial step in such type of variational methods is to ensure the compactness conditions  stated in the following lemma.
\begin{lemma}
Under the assumptions of Theorem $\ref{th3.1}$, the functional $\mathcal{J}_{\lambda}$ verifies the Palais-Smale condition.
\end{lemma}
\noindent{\bf{\textit{Proof}.}}
Let $\{w_n\}$ be a (PS)-sequence of $\mathcal{J}_{\lambda}$ on $X_0$, which could be expressed mathematically as follows 
 \begin{equation}
      \mathcal{J}_{\lambda}\left(w_{n}\right) \rightarrow c>0 \quad \text{and} \quad \mathcal{J}_{\lambda}^{\prime}\left(w_{n}\right) \rightarrow 0 \text { in } X_0^{*}.
      \label{Jbound}
      \end{equation}
Using Lemma \ref{l2.1},  we have $$
 \begin{aligned}
\mathcal{J}_{\lambda} \left(w_{n}\right) &= \sum_{i=1}^N \int_{Q} \tfrac{1}{p_i(x, y)} \frac{|w_n(x)-w_n(y)|^{p_i(x, y)}}{|x-y|^{s p_i(x, y)+N}} d x d y \\
& ~~+  \int_{\Omega}
 \tfrac{1}{p_M(x)} |w_n(x)|^{p_M(x)} dx - \lambda \int_{\Omega} \frac{1}{r(x)}  |w_n(x)|^{r(x)} dx \\
& \geqslant \frac{1}{P^+_+} \sum_{i=1}^N \int_{Q} \frac{|w_n(x)-w_n(y)|^{p_i(x, y)}}{|x-y|^{s p_i(x, y)+N}} d x dy + \frac{1}{p^+_M} \int_{\Omega}
 |w_n(x)|^{p_M(x)} dx \\
 &-   \frac{\lambda}{r^-} \int_{\Omega} |w_n(x)|^{r(x)} dx \\
& \geqslant \min \left\lbrace  \tfrac{1}{P^+_+},\tfrac{1}{p^+_M}\right\rbrace  \left( \rho_{ \vec{p}(x,y)}^0 (w_n) + \rho_{p_M(x)} (w_n)\right) - \tfrac{\lambda}{r^-} \rho_{r(x)} (w_n) \\
& \geqslant  \min \left\lbrace  \tfrac{1}{P^+_+},\tfrac{1}{p^+_M}\right\rbrace \left( \sum_{i=1}^N [w_n]_{s,p_i(x, y)}^{p_i^-} + \|w_n\|_{p_M(x)}^{p_M^-} \right) - \tfrac{\lambda}{r^-}  \|w_n\|_{r(x)}^{r^+} \\
& \geqslant \min \left\lbrace  \tfrac{1}{P^+_+},\tfrac{1}{p^+_M}\right\rbrace \left( \frac{[w_n]_{s,\vec{p}(x, y)}^{P_-^-}}{N^{P^-_- -1}}  + \|w_n\|_{p_M(x)}^{p_M^-} \right) - \tfrac{\lambda}{r^-}  \|w_n\|_{r(x)}^{r^+}\\
& \geqslant \min \left\lbrace  \tfrac{1}{P^+_+},\tfrac{1}{p^+_M}\right\rbrace \tfrac{1}{N^{P^-_- -1} 2^{\tilde{P}-1}} \|w_n\|_{X}^{\tilde{P}} - \frac{\lambda}{r^-}  \|w_n\|_{r(x)}^{r^+},
 \end{aligned}
 $$
 where $\tilde{P}= \min\left\lbrace  P^-_-, P_M^-\right\rbrace$, and by using Theorem \ref{th2}, we get that 
 \begin{equation}
  \mathcal{J}_{\lambda} \left(w_{n}\right) \geqslant \min\left\lbrace \tfrac{1}{P^+_+},\tfrac{1}{p^+_M}\right\rbrace  \frac{1}{N^{P^-_- -1} 2^{\tilde{P}-1}} [w_n]_{s, \vec{p}(x, y)}^{\tilde{P}}-  
 \tfrac{\lambda C^{r^+}}{r^-} [w_n]_{s, \vec{p}(x, y)}^{r^+}.
 \label{minoration}  \end{equation}
  Since $\mathcal{J}_{\lambda}$ is bounded from \eqref{Jbound}, thus $\{w_n\}$ is bounded on $X_0$. In the view of the reflexivity of $X_0$, we can find a subsequence of $\{w_n\}$ still denoted in the same way such that $w_n \rightharpoonup w$ in $X_0$. \\
  From condition \eqref{condition1}, $X_0$ is compactly embedded in $L^{r(x)}(\Omega)$ (resp. in $L^{p_{\scriptsize M}(x)}(\Omega)$), therefore $w_n \longrightarrow w$ in $L^{r(x)}(\Omega)$ (resp. in $L^{p_{\scriptsize M}(x)}(\Omega)$), using H\"{o}lder's inequality we acquire 
 \begin{equation}
 \int_{\Omega}\left|w_n\right|^{r(x)-2} w_n\left(w_n-w\right) d x \leqslant 2\left\||w_n|^{r(x)-1}\right\|_{L^{r'(x)}(\Omega)}\left\|w_n-w\right\|_{L^{r(x)}(\Omega)},
 \label{3.4}
 \end{equation}
 Similarly, 
 \begin{equation}
  \int_{\Omega}\left|w_n\right|^{p_M(x)-2} w_n\left(w_n-w\right) d x \leqslant 2\left\||w_n|^{p_M(x)-1}\right\|_{L^{p_M'(x)}(\Omega)}\left\|w_n-w\right\|_{L^{p_M(x)}(\Omega)}.
  \label{4.5}
  \end{equation}
  Hence 
  $$\lim_{n\to \infty}\int_{\Omega}\left|w_n\right|^{r(x)-2} w_n\left(w_n-w\right) d x=\lim_{n\to \infty}\int_{\Omega}\left|w_n\right|^{p_{\scriptsize M}(x)-2} w_n\left(w_n-w\right) d x=0.$$
  Owing to \eqref{Jbound}, and the weak convergence of ${w_n}$ we have 
  $$\left\langle \mathcal{J}^{\prime}_{\lambda}\left(w_{n}\right), w_{n}-w\right\rangle \longrightarrow 0 \quad \text{as } n \longrightarrow +\infty.$$
  Which means that 
{\small   $$     
  \begin{aligned}
\sum_{i=1}^N & \int_{Q} \frac{|w_n(x)-w_n(y)|^{p_i(x, y)-2}(w_n(x)-w_n(y))((w_n-w)(x)-(w_n-w)(y))}{|x-y|^{s p_i(x, y)+N}} d x d y \\
&+\int_{\Omega}|w_n(x)|^{p_M(x)-2} w_n(x) (w_n-w)(x) d x-\lambda \int_{\Omega}|w_n(x)|^{r(x)-2} w_n(x) ( w_n-w)(x) d x  \underset{n\longrightarrow +\infty}{\longrightarrow 0}.
  \end{aligned}
  $$}
From \eqref{3.4}, and \eqref{4.5}, it follows that 
{\small $$
\sum_{i=1}^N \int_{Q} \frac{|w_n(x)-w_n(y)|^{p_i(x, y)-2}(w_n(x)-w_n(y))((w_n-w)(x)-(w_n-w)(y))}{|x-y|^{s p_i(x, y)+N}} d x d y 
 \underset{n\longrightarrow +\infty}{\longrightarrow 0}.
$$}
On the other hand, it's clear that 
{\small $$
\sum_{i=1}^N \int_{Q} \frac{|w(x)-w(y)|^{p_i(x, y)-2}(w(x)-w(y))((w_n-w)(x)-(w_n-w)(y))}{|x-y|^{s p_i(x, y)+N}} d x d y 
 \underset{n\longrightarrow +\infty}{\longrightarrow 0}.
$$}
The two previous results show that  $\left\langle \mathcal{L}_{\vec{p}(x,y)} \left( w_n\right)- \mathcal{L}_{\vec{p}(x,y)} \left( w_n\right), w_n-w\right\rangle \underset{n\longrightarrow +\infty}{\longrightarrow 0} $.
The fact that $\mathcal{L}_{\vec{p}(.,.)}$ is of type $(S_+)$ implies the strong convergence of $w_n$ to $w$ in $X_0.$ We deduce that $\mathcal{J}_{\lambda}$ satisfies the Palais-Smale condition.\\

Now, let us verify that the functional $\mathcal{J}_{\lambda}$
  satisfies the mountain pass geometry.
\begin{lemma}
 Assume that the hypothesis of Theorem $\ref{existence}$ are fulfilled. Then, there exists $\lambda_*>0$ such that for any $\lambda \in\left(0, \lambda_*\right)$, there are $\delta, \vartheta >0$, such that $\mathcal{J}_\lambda(w) \geqslant \vartheta>0$ for any $w \in X_0 $ with $\|w\|_{X_0}=\delta$.
 \label{lemma2}
\end{lemma}
\noindent{\bf{\textit{Proof}.}}
From inequality \eqref{minoration}, we have:
$$
\mathcal{J}_{\lambda} \left(w\right) \geqslant \tilde{C} [w]_{s, \vec{p}(x, y)}^{\tilde{P}}-  
 \frac{\lambda C^{r^+}}{r^-} [w]_{s, \vec{p}(x, y)}^{r^+},~ \forall w \in X_0,
$$
where $\tilde{C} = \min\left\lbrace \tfrac{1}{P^+_+},\tfrac{1}{p^+_M}\right\rbrace  \tfrac{1}{N^{P^-_- -1} 2^{\tilde{P}-1}}$. Therefore, for any $w \in X_0$ with $\|w\|_{X_0}=\delta$, we obtain
$$
\mathcal{J}_{\lambda} \left(w\right) \geqslant \tilde{C} \delta^{\tilde{P}}-  
 \tfrac{\lambda C^{r^+}}{r^-}\delta^{r^+} = \delta^{r^+}\left( \tilde{C} \delta^{\tilde{P}-r^+}- \tfrac{\lambda C^{r^+}}{r^-}\right). 
$$
By defining $\lambda_*$ in this way $ \lambda_*=  \frac{\bar{C}\delta^{\bar{P}-r^+}r^-}{2C^{r^+}}  $, we get as a result that  for all $ \lambda \in (0,\lambda_*),$ for any  $w \in X_0$ with $\|w\|_{X_0}=\delta$ such that 
$$ \mathcal{J}_{\lambda}(w) \geqslant \vartheta >0  \quad \text{with } \vartheta= \frac{\delta}{2}. $$
\begin{lemma}
 Assume that the hypothesis of Theorem $\ref{existence}$ are fulfilled. Then, there exists $\omega \in X_0$ such that $\omega \geqslant 0, \omega \neq 0$ and $\mathcal{J}_\lambda(t \omega)<0$ for any $t>0$ small enough.
 \label{lemma3}
\end{lemma}
\noindent{\bf{\textit{Proof}.}}
Since $r^- <P^-_-$, then there exists $ \varepsilon >0$ such that $r^- + \varepsilon \leqslant P^-_-$. Moreover, the continuity of $r$ implies the existence of an open set $O$ of $\Omega$ such that \begin{equation}
    |r(x)-r^-| \leqslant \varepsilon \quad\text{for all}~~ x \in O. \end{equation}
Consider $\omega$ as an infinitely differentiable function that is compactly supported such that $\overline{O} \subset supp (\omega)$, $\omega(x)=1 ~~ \forall x \in \overline{O}$, and $0 \leqslant \omega \leqslant 1 ~~ \text{in } \Omega$. Therefore, for $t>0$ small enough we have:
$$
\begin{aligned}
\mathcal{J}_{\lambda}(t \omega)&=  \sum_{i=1}^N \int_{Q} \frac{1}{p_i(x, y)} \frac{|t \omega(x)-t \omega(y)|^{p_i(x, y)}}{|x-y|^{s p_i(x, y)+N}} d x d y +  \int_{\Omega}
 \frac{1}{p_M(x)} |t \omega(x)|^{p_M(x)} dx \\
& - \lambda \int_{\Omega} \frac{1}{r(x)}  |t \omega(x)|^{r(x)} dx \quad \\
& \leqslant  \frac{t^{P^-_-}}{P^-_-} \sum_{i=1}^N \int_{Q} \frac{|\omega(x)- \omega(y)|^{p_i(x, y)}}{|x-y|^{s p_i(x, y)+N}} d x d y + \frac{t^{p^-_{M}}}{p^-_{M}} \int_{\Omega}
  |\omega(x)|^{p_M(x)} dx \\
  & - \frac{\lambda }{r^+}  \int_{\Omega}   t^{r(x)}|\omega(x)|^{r(x)} dx, \\
  & \leqslant \frac{t^{P^-_-}}{P^-_-} \rho_{s,\vec{p}(x,y)}(\omega) -\frac{\lambda t^{r^- + \varepsilon}}{r^+} \int_{O}  |\omega(x)|^{r(x)} dx, \\
  & = t^{r^-+\varepsilon} \left( \frac{t^{P^-_- - r^- - \varepsilon}}{P^-_-} \rho_{s,\vec{p}(x,y)}(\omega)- \frac{\lambda }{r^+} \int_{O}  |\omega(x)|^{r(x)} dx\right).
\end{aligned}
$$
For $t< \theta^{\frac{1}{P^-_- - r^- - \varepsilon}}$, we obtain that $\mathcal{J}_{\lambda}(t \omega)<0$, where $$
0<\theta < \min \left\{ 1, \frac{\lambda P^-_- \int_{O}  |\omega(x)|^{r(x)} dx}{r^+ \rho_{s,\vec{p}(x,y)} (\omega)}  \right\},$$
Next, we point-out that $\rho_{s,\vec{p}(x,y)} (\omega)>0$ (this fact implies that $w\neq0$).
Indeed, from this inclusion between sets $O \subset \operatorname{supp}(\omega) \subset \Omega$ and $0 \leqslant \omega \leqslant 1$ in $\Omega$, we acquire
$$
0<\int_{0}|\omega|^{r(x)} \mathrm{d} x \leqslant \int_{\Omega}|\omega|^{r(x)} \mathrm{d} x \leqslant \int_{\Omega}|\omega|^{r^{-}} \mathrm{d} x .
$$
Condition \eqref{condition1} implies the existence of a constante $C$ such that 
$$
\|\omega\|_{L^{r^{-}(\Omega)}} \leqslant C\|\omega\|_{X_0} .
$$
The two previous inequalities inferredes to 
$$ \rho_{s,\vec{p}(x,y)} (\omega) >0. $$
The proof of Lemma \ref{lemma3} is achieved.\\
\noindent Now, we are in the position to prove our existence result.\\

\noindent{\bf{\textit{Proof of Theorem} \ref{th3.1}.}} From Lemma \ref{lemma2} we have: 
\begin{equation}
 \inf \mathcal{J}_\lambda>0   \qquad ~~~~ \text{  in  } \partial \mathcal{B}_m(0), 
\end{equation}
and from Lemma \ref{lemma3}, we get  \begin{equation}
 \mathcal{J}_{\lambda} \left(w\right) \geqslant \tilde{C} [w]_{s, \vec{p}(x, y)}^{\tilde{P}}-  
 \frac{\lambda C^{r^+}}{r^-} [w]_{s, \vec{p}(x, y)}^{r^+},~ \forall w \in   \mathcal{B}_m(0), 
\end{equation}
which implies that $-\infty < c = \inf _{\overline{\mathcal{B}_m(0)} } \mathcal{J}_\lambda <0$. \\
Consider $\varepsilon> 0$ such that $$0 < \varepsilon < \inf _{\partial \mathcal{B}_m(0) } \mathcal{J}_\lambda-\inf _{\mathcal{B}_m(0)} \mathcal{J}_\lambda.$$
Then we apply the Ekelnad's variational principle to the functional $ \mathcal{J}_{\lambda} : \overline{\mathcal{B}_m(0)} \longrightarrow \mathbb{R}$, we get $w_{\varepsilon} \in \overline{\mathcal{B}_m(0)}$ such that 
$$
\left\{\begin{array}{l}
\mathcal{J}_\lambda\left(w_{\varepsilon}\right) \leq  \inf _{\overline{\mathcal{B}_m(0)}} \mathcal{J}_\lambda+\varepsilon, \\
\mathcal{J}_\lambda\left(w_{\varepsilon}\right)<\mathcal{J}_\lambda(w)+\varepsilon\left\|w-w_{\varepsilon}\right\|_{X_0}, \quad \text { for } w \neq w_{\varepsilon}, w_{\varepsilon} \in \mathcal{B}_m(0)  .
\end{array}\right.
$$
Which acquires that 
$$ \mathcal{J}_\lambda\left(w_{\varepsilon}\right) \leq \inf _{\mathcal{B}_m(0)} \mathcal{J}_\lambda+\varepsilon \leq \inf _{\mathcal{B}_m(0)} \mathcal{J}_\lambda+\varepsilon<\inf _{\partial \mathcal{B}_m(0) } \mathcal{J}_\lambda. $$
In following, we cosndier the mapping 
$$
\begin{aligned}
\mathcal{I}_\lambda^{\varepsilon}: \overline{\mathcal{B}_m(0)} & \longrightarrow \mathbb{R} \\
w & \longrightarrow \mathcal{J}_\lambda(w)+\varepsilon\left\|w-w_{\varepsilon}\right\|_{X_0} .
\end{aligned}
$$
It is clear that $\mathcal{I}_\lambda^{\varepsilon} (w_\varepsilon) = \mathcal{J}_\lambda(w_\varepsilon) <  \mathcal{J}_\lambda(w) $, for $w \neq w_\varepsilon$. Hence, $w_\varepsilon$ is a minimum of $ \mathcal{I}_\lambda^{ \varepsilon}$. Which implies that for any $t \in (0,1)$ and $v \in \mathcal{B}_{m}$, we have 
$$ \frac{\mathcal{I}_\lambda^{\varepsilon}\left(w_{\varepsilon}+t v \right)-\mathcal{I}_\lambda^{\varepsilon}\left(w_{\varepsilon}\right)}{t} \geqslant 0 \Rightarrow  \frac{\mathcal{J}_\lambda\left(w_{\varepsilon}+t v\right)-\mathcal{J}_\lambda\left(w_{\varepsilon}\right)}{t}+\varepsilon\|v\|_{X_0} \geqslant 0,  $$
as $t$ tends to $0^+$, we acquire $ <\mathcal{J}_\lambda^{\prime}\left(w_{\varepsilon}\right), v> + \varepsilon \|v\|_{X_0}   \geqslant 0 .$ \\
As a consequence, $$ \|\mathcal{J}_\lambda^{\prime} \left(w_{\varepsilon}\right)\|_{X_0^*} \leqslant \varepsilon, $$ 
we infer that there exists a sequence $w_n \in \mathcal{B}_m$  such that 
$$  \mathcal{J}_{\lambda}\left(w_{n}\right) \rightarrow c>0 \quad \text{and} \quad \mathcal{J}_{\lambda}^{\prime}\left(w_{n}\right) \rightarrow 0 .  $$
$ \mathcal{J}_{\lambda}$ verifies the PS condition, then we can find $w_0 \in X_0 $ such that $w_n \longrightarrow w_0$ as $n \longrightarrow +\infty$. as a result, $\mathcal{J}_{\lambda} (w_0)=c$ and $ \mathcal{J}_{\lambda}^{\prime} (w_0)=0$. In conculsion, $w_0$ is a critical point  of $\mathcal{J}_{\lambda}$. Which is summed up by the achievement of the proof.\\\vspace*{0.5cm}

\noindent\textbf{{\large Declarations}} :\\

\noindent\textbf{Ethical Approval :} Not applicable.\\

\noindent\textbf{Competing interests :} The authors declare that there is no conflict of interest.\\

\noindent\textbf{Authors' contributions :} The authors contributed equally to this work.\\

\noindent\textbf{Funding :} Not applicable.\\

\noindent\textbf{Availability of data and materials :} Not applicable.


\begin{thebibliography}{99}
\bibitem{ambrosetti} A. Ambrosetti and P. H. Rabinowitz, Dual variational methods in critical point theory and applications, J. Funct. Anal., 14 (1973), 349-341.
\bibitem{ayazoglu} R. Ayazoglu, Y. Saraç, S. \c{S}. \c{S}ener, and G. Alisoy, Existence and multiplicity of solutions for a Schr\"{o}dinger–Kirchhoff type equation involving the fractional $p(.,.)$-Laplacian operator in $\mathbb{R^N}$, Collectanea Mathematica, vol. 72, pp. 129-156 (2020).
\bibitem{azroul2} E. Azroul, A. Benkirane, {M. Shimi},  and  M. Srati,  \textit{On a class of fractional  $p\left( x\right)$-Kirchhoff type problems}, J. Applicable Analysis, vol. 100 (2021), no. 2, pp. 383-402,  \url{https://doi.org/10.1080/00036811.2019.1603372}.
  \bibitem{azroul1} E. Azroul,  A. Benkirane, and  M. Shimi,  \textit{Eigenvalue problems involving the fractional  $p(x)$-Laplacian operator}, J. Advances in Operator Theory, vol. 4 (2019), no. 2, pp. 539--555. \url{https://doi.org/10.15352/aot.1809-1420}
\bibitem{azroul3} E. Azroul,  A. Benkirane, and  {M. Shimi}, \textit{General fractional Sobolev space with variable exponent and applications to nonlocal problems}, Advances in Operator Theory, vol. 5 (2020), no 4, pp. 1512-1540. \url{https://doi.org/10.1007/s43036-020-00062-w}. 
 \bibitem{azroul4} E. Azroul,  A. Benkirane,  N. T. Chung, and  {M. Shimi}, \textit{Existence results for anisotropic fractional $\left( p_1\left( x,.\right) ,p_2\left( x,.\right)\right)  $-Kirchhoff type problems}, Journal of Applied Analysis and Computation, vol. 11 (2021), no. 5, pp. 1-24, \url{doi: 10.11948/20200394}.
  \bibitem{bahrouni} A. Bahrouni and V. R\u{a}dulescu, On a new fractional Sobolev space and applications to
  nonlocal variational problems with variable exponent, Discrete Contin. Dyn. Syst. Ser. S 11 (2018), no. 3, 379-389.
   \bibitem{bedahmane} M. Bendahmane, M. Chrif and S. El Manouni, An approximation result in
   generalized anisotropic Sobolev spaces and application, Z. Anal. Anwend., 30(3)
   (2011), 341-353.
 \bibitem{benedek} A. Benedek and R. Panzone, The space $L^p$, with mixed norm, Duke Math. J. 28 (1961), pp. 301-324.
  \bibitem{brezis} H. Brezis, Functional Analysis, Sobolev Spaces and Partial Diffrential Equations, Universitext,  Springer, New York, 2011.
 \bibitem{chen}W. Chen and S. Deng,  The Nehari manifold for a fractional p-Laplacian system involving concave-convex nonlinearities, Nonlinear Anal, vol. 27, pp. 80–92, 2016.
 
 \bibitem{diening} L. Diening, P. Harjulehto, P. Hästö, M. Ruzicka, Lebesgue and Sobolev Spaces with Variable Exponents, vol. 2017,
 Springer, Heidelberg, 2011.
 \bibitem{dinezza} E. Di Nezza, G; Palatucci, and E. Valdinoci , \textit{Hitchhiker's guide to the fractional Sobolev spaces}, Bulletin des sciences math{\'e}matiques, vol. 2 (2011), no. 3.
  \bibitem{engheta}
   N.Engheta, On the role of fractional calculus in electromagnetic theory,
 IEEE Antennas Propagat. Mag., vol. 39, no. 4, pp. 35–46, 1997
  \bibitem{fan1} X. Fan, Anisotropic variable exponent Sobolev spaces and
  $\vec{p}(.)$-Laplacian equations, Complex Var. Elliptic Equ. 56 (2011) 623–642.
   \bibitem{zhang} X. L. Fan and Q. H. Zhang, Existence of solutions for $p(x)$-Laplacian Dirichlet problem,
   Nonlinear Anal. 52 (2003), 1843-1852.
        \bibitem{Fan} X. L. Fan, D. Zhao, On the spaces $L^{p(x)}(\Omega)$ and $W^{m,p(x)}(\Omega)$ , J. Math. Anal. Appl. 263
        (2001), 424-446. 
 \bibitem{ho} K. Ho, Y. H. Kim, A-priori bounds and multiplicity of solutions for nonlinear elliptic problems
 involving the fractional $p(.)$-Laplacian, Nonlinear Anal., 188 (2019), 179--201.
 \bibitem{hormander}
 L. H\"{o}rmander, Estimates for translation invariant operators in $L^p$
 spaces, Acta Math. 104
 (1960), 93-140.
  \bibitem{igari}  S. Igari, Interpolation of operators in Lebesgue spaces with mixed norm and its applications
 to Fourier analysis, Tohoku Math. J. (2) 38 (1986), pp. 469-490.
 \bibitem{kaufmann}U. Kaufmann, J.D.  Rossi, and R.E. Vidal, \textit{Fractional Sobolev spaces with variable exponents and fractional $p (x)$-Laplacians}, Electronic Journal of Qualitative Theory of Differential Equations, vol. (2017), no. 76, pp. 1--10. 
  \bibitem{Kovacik}O. Kov{\'a}cik and J. R{\'a}kosn{\i}k. On Spaces $L^{p(x)}(\Omega)$ and $W^{m,p(x)}(\Omega)$, Czechoslovak Math, 41(116):592-618, 1991.
 \bibitem{kolodi}] S.N. Kruzhkov, I.M. Kolodii, On the theory of embedding of anisotropic Sobolev spaces, Russian Math. Surveys 38 (1983) 188–189.
\bibitem{Liu} H. K. Liu, Y. Q. Fu, On the variable exponential fractional Sobolev space $W^{s(.),p(.)}$, AIMS Math., 6
 (2020), 6261--6276.
 \bibitem{madych}W. R. Madych, Lipschitz spaces and mixed Lebesgue spaces, Proc. Amer. Math. Soc. 85
 (1982), pp. 213-218. 
  \bibitem{mihalescu}
  M. Mih\u{a}ilescu, P. Pucci and V. D. R\u{a}dulescu, Nonhomogeneous boundary value problems in anisotropic Sobolev spaces, C. R. Acad. Sci. Paris, Ser. I 345 (2007), 561-566.
   \bibitem{mihalescu2} M. Mih\u{a}ilescu and V. D. R\u{a}dulescu, On a nonhomogeneous quasilinear eigenvalue problem
   in Sobolev spaces with variable exponent, Proc. Amer. Math. Soc. 135 (2007), no. 9, 2929-2937.
   \bibitem{Musielak} J. Musielak, Orlicz spaces and Modular spaces, Lecture Notes in Mathematics, vol. 1034. Springer-Verlag Berlin, 1983.
   
   \bibitem{Natanson}I.P. Natanson, Theory of Functions of a Real Variable. Nauka, Moscow (1950).
  \bibitem{Oldham}K. B. Oldham, Fractional differential equations in electrochemistry, Advances in Engineering Software, vol. 41, no. 1, pp. 9-12, 2010.
  \bibitem{servadei}R. Servadei and E. Valdinoci, Lewy-Stampacchia type estimates for variational inequalities driven by
  nonlocal operators, Rev. Mat. Iberoam. vol. 29, no. 3, pp. 1091–1126, 2013. 
 \bibitem{simon} J. Simon, Régularité de la solution d'une équation non linéaire dans $\mathbb{R}^N$, Journées d'Analyse Non Linéaire (Proc. Conf., Besançon, 1977), Lecture Notes in Math., 665, Springer, Berlin, 1978, pp. 205-227. 
  \bibitem{sjodin}T. Sj\"{o}din, Weighted norm inequalities for Riesz potentials and fractional maximal functions
  in mixed norm Lebesgue spaces, Studia Math. 97 (1991), pp.  239-244.
   \bibitem{Sobolev} S.L. Sobolev, On a theorem of functional analysis, Mat. Sb. 4 (46) 1938, 39--68 (translated into English
   in 1963).
  \bibitem{troisi} M. Troisi, Teoremi di inclusione per spazi di Sobolev non isotropi. (Italian) Ricerche Mat. 18
  (1969), pp. 3-24. 
   \bibitem{Xu} C. Xu and W. Sun,  An embedding theorem for anisotropic fractional Sobolev spaces. Banach J. Math. Anal. 15, 63 (2021). https://doi.org/10.1007/s43037-021-00146-6.
  \bibitem{yang} Q. Yang, D. Chen, T. Zhao, and Y. Chen. Fractional calculus in image processing: a review. Frac. Calc. Appl. Anal.,
  19(5):1222–1249, 2016.  
   \bibitem{Zuo}  J. Zuo, T. An, A. Fiscella,  A critical Kirchhoff-type problem driven by a $p(.)$-fractional
 Laplace operator with variable $s(.)$-order. Math Meth Appl Sci. 2020;1--15. \url{https://doi.org/10.1002/mma.6813}
\end{thebibliography}
\end{document}